\documentclass[12pt]{article}
\usepackage{latexsym,amssymb,amsmath}
\begin{document}

\baselineskip=20pt

\newcommand{\la}{\langle}
\newcommand{\ra}{\rangle}
\newcommand{\psp}{\vspace{0.4cm}}
\newcommand{\pse}{\vspace{0.2cm}}
\newcommand{\ptl}{\partial}
\newcommand{\dlt}{\delta}
\newcommand{\sgm}{\sigma}
\newcommand{\al}{\alpha}
\newcommand{\be}{\beta}
\newcommand{\G}{\Gamma}
\newcommand{\gm}{\gamma}
\newcommand{\vs}{\varsigma}
\newcommand{\Lmd}{\Lambda}
\newcommand{\lmd}{\lambda}
\newcommand{\td}{\tilde}
\newcommand{\vf}{\varphi}
\newcommand{\yt}{Y^{\nu}}
\newcommand{\wt}{\mbox{wt}\:}
\newcommand{\rd}{\mbox{Res}}
\newcommand{\ad}{\mbox{ad}}
\newcommand{\stl}{\stackrel}
\newcommand{\ol}{\overline}
\newcommand{\ul}{\underline}
\newcommand{\es}{\epsilon}
\newcommand{\dmd}{\diamond}
\newcommand{\clt}{\clubsuit}
\newcommand{\vt}{\vartheta}
\newcommand{\ves}{\varepsilon}
\newcommand{\dg}{\dagger}
\newcommand{\tr}{\mbox{Tr}}
\newcommand{\ga}{{\cal G}({\cal A})}
\newcommand{\hga}{\hat{\cal G}({\cal A})}
\newcommand{\Edo}{\mbox{End}\:}
\newcommand{\for}{\mbox{for}}
\newcommand{\kn}{\mbox{ker}}
\newcommand{\Dlt}{\Delta}
\newcommand{\rad}{\mbox{Rad}}
\newcommand{\rta}{\rightarrow}

\begin{center}{\LARGE \bf Gel'fand-Dorfman Bialgebras}\end{center}
\vspace{0.2cm}

\begin{center}{\large Xiaoping Xu}\end{center}

(A talk given in the workshop on Algebra and Discrete Mathematics at The Chinese University of Hong Kong, 27 - 31 March 2000)  

\section{Introduction}

A Gel'fand-Dorfman bialgebra is a vector space with a Lie algebra structure and a Novikov algebra structure, satisfying a certain compatibility condition. This bialgebraic structure corresponds to a certain Hamiltonian pairs in the Gel'fand-Dikii-Dorfman's theory of Hamiltonian operators (cf. [GDi], [GDo]). In this talk, I will give a survey on the study of Gel'fand-Dorfman bialgebras. First let me give a more technical initial introduction to the bialgebras through what I call  ``Lie algebra with one-variable
structure.''

Let  $\Bbb{F}$ be an arbitrary field. Denote by $\Bbb{Z}$ the ring of integers and by $\Bbb{N}$ the additive semi-group of nonnegative integers. All the vector spaces are assumed over $\Bbb{F}$. Let $V$ be a vector space and let $t$ be an indeterminant. Form a tensor
$$\hat{V}=V\otimes_{\Bbb{F}}\Bbb{F}[t,t^{-1}].\eqno(1.1)$$
We denote
$$u(z)=\sum_{n\Bbb{Z}}(u\otimes t^n)z^{-n-1}\qquad\for\;\;u\in V,\eqno(1.2)$$
where $z$ is a formal variable. A {\it Lie algebra with one-variable structure} is a vector space $\hat{V}$ with the Lie bracket of the form
$$[u(z_1),v(z_2)]=\sum_{i=0}^m\sum_{j=0}^nz_2^{-1}\ptl_{z_1}^i\dlt\left({z_1\over z_2}\right)\ptl_{z_2}^jw_{ij}(z_2)\eqno(1.3)$$
for $u,v\in V$, where $m,\;n$ are nonnegative integers depending on $u,v$, and $w_{ij}\in V$. Here we have used the notation
 $$\dlt(z)=\sum_{k\in\Bbb{Z}}z^k.\eqno(1.4)$$
Moreover, each $w_{ij}$ depends on $u$ and $v$ bilinearly by (1.2) and (1.3). Above definition is not in general form. We aim at giving the audience a simple rough picture. According to Wightman's axioms of quantum field theory,
the algebraic content of two-dimensional quantum field theory is a certain representation theory of Lie algebras with one-variable structure (e.g., cf. [K]),
where intertwining operators among irreducible modules, ``partition functions'' (characters in algebraic terms) and ``correlation functions'' related to (1.3) etc.  play important roles.

Suppose that $m=n=0$ in (1.3). We denote
$$ w_{00}=[u,v].\eqno(1.5)$$
Then $(V,[u,v])$ forms a Lie algebra and $\hat{V}$ is the corresponding {\it loop algebra} with the Lie bracket
$$[u\otimes t^j,v\otimes t^k]=[u,v]\otimes t^{j+k}\qquad\for\;\;u,v\in V,\;j,k\in\Bbb{Z}.\eqno(1.6)$$

If $V=\Bbb{F}e$ is one-dimensional and
$$[e(z_1),e(z_2)]=z_2^{-1}\dlt\left({z_1\over z_2}\right)\ptl_{z_2}e(z_2)-2z_2^{-1}\ptl_{z_1}\dlt\left({z_1\over z_2}\right)e(z_2),\eqno(1.7)$$
Then $\hat{V}$ is the {\it centerless Virasoro algebra} (or a rank-one Witt algebra) with the Lie bracket
$$[e\otimes t^j,e\otimes t^k]=(j-k)e\otimes t^{j+k}\qquad\for\;\;j,k\in\Bbb{Z}.\eqno(1.8)$$

A natural generalization of the above one-dimensional case is
$$[u(z_1),v(z_2)]=z_2^{-1}\dlt\left({z_1\over z_2}\right)\ptl_{z_2}w_{01}(z_2)+z_2^{-1}\ptl_{z_1}\dlt\left({z_1\over z_2}\right)w_{10}(z_2)\eqno(1.9)$$
for $u,v\in V$. Denote
$$w_{01}=v\circ u.\eqno(1.10)$$
Then $(V,\circ)$ is an algebra satisfying
$$(u\circ v)\circ w=(u\circ w)\circ v,\eqno(1.11)$$
$$(u\circ v)\circ w-u\circ (v\circ w)=(v\circ u)\circ w-v\circ (u\circ w)\eqno(1.12)$$
for $u,v,w\in V.$ The above algebra $(V,\circ)$ appeared in Gel'fand and Dorfman's work [GDo], corresponding to certain Hamiltonian operators. Moreover, it also appeared in Balinskii and Novikov's work [BN]  as the local structures of certain Poisson brackets of hydrodynamic type. This structure was first abstractly studied by Zel'manov [Z], Filippov [F] and was named as ``Novikov algebra" by Osborn [O1]. 

Note that (1.12) is the axiom of left-symmetric algebra. Left-symmetric algebras play fundamental roles in the theory of affine manifolds (cf. [A], [FD]). Novikov algebras are Left-symmetric algebras whose right multiplication operators  are mutually commutative (cf. (1.11)). Furthermore,
$$w_{10}=-(u\circ v+v\circ u).\eqno(1.13)$$
The Lie bracket on $\hat{V}$ is
$$[u\otimes t^j,v\otimes t^k]=ju\circ v\otimes t^{j+k-1}-kv\circ u\otimes t^{j+k-1}
\qquad\for\;\;u,v\in V,\;j,k\in\Bbb{Z}.\eqno(1.14)$$

Consider the mixed case
$$[u(z_1),v(z_2)]=z_2^{-1}\dlt\left({z_1\over z_2}\right)[w_{00}(z_2)+\ptl_{z_2}w_{01}(z_2)]+z_2^{-1}\ptl_{z_1}\dlt\left({z_1\over z_2}\right)w_{10}(z_2)\eqno(1.15)$$
for $u,v\in V$. Denote
$$w_{00}=[v,u],\qquad w_{01}=v\circ u.\eqno(1.16)$$
Then $(V,[\cdot,\cdot])$ forms a Lie algebra, $(V,\circ)$ forms a Novikov algebra and the following compatibility condition holds
$$[w\circ u,v]-[w\circ v,u]+[w,u]\circ v-[w,v]\circ u-w\circ [u,v]=0\eqno(1.17)$$
for $u,v,w\in V$.

I called the triple $(V,[\cdot,\cdot],\circ)$ a {\it Gel'fand-Dorfman bialgebra} (cf. [X4]), which corresponds to a certain Hamiltonian pair in Gel'fand and Dorfman's work [GDo]. Moreover, (1.13) holds. The Lie bracket on $\hat{V}$ is
$$[u\otimes t^j,v\otimes t^k]=[v,u]\otimes t^{j+k}+ju\circ v\otimes t^{j+k-1}-kv\circ u\otimes t^{j+k-1}
\eqno(1.18)$$
 for $u,v\in V,\;j,k\in\Bbb{Z}$. 

For convenience, we shall use the following notation of index
$$\ol{m,n}=\{m,m+1,...,n\}\eqno(1.19)$$
for $m,n\in\Bbb{N}$ such that $m\leq n$.

The article is organized as follows. In Section 2, I will talk about  structures of simple Novikov algebras and their irreducible representations. In Section 3, I will present some general constructions of Gel'fand-Dorfman bialgebras.
In Section 4, I will give a few classification results on the bialgebras. In Section 5,  examples of application to simple ``cubic conformal algebras'' and ``quartic conformal algebras'' will be given.

\section{Novikov Algebras}

Let us first give a construction of Novikov algebras. Let $({\cal A},\cdot)$ be a commutative associative algebra. Then $({\cal A},\cdot)$ forms a Novikov algebra, which is not so interesting. Take a derivation $\ptl$ of $({\cal A},\cdot)$ and $\xi\in {\cal A}$. We define
$$u\circ_{\xi}v=u\ptl(v)+\xi uv\qquad\for\;\;u,v\in{\cal A}.\eqno(2.1)$$
Then $({\cal A},\circ_{\xi})$ forms a Novikov algebra. The above construction was found by S. Gel'fand when $\xi=0$ (cf. [GDo]), by Fillip [F] when $\xi\in\Bbb{F}$ and by me [X2] in general case. 
\psp

{\bf Theorem 2.1 (Zel'manov [Z])}. {\it  Any finite-dimensional simple Novikov algebra over an algebraically closed field} $\Bbb{F}$ {\it with characteristic} $0$ {\it isomorphic to} $(\Bbb{F},\cdot)$.
\psp

Based on Osborn's work [O1], I obtained the following result:
\psp 

{\bf Theorem 2.2 (Xu, [X1])}. {\it A finite-dimensional simple Novikov algebra over  an algebraically closed field} $\Bbb{F}$ {\it with characteristic} $p>2$  {\it is either one-dimensional with a base element} $e$ {\it such that} $e\circ e=e$ {\it or  has dimension} $p^k$ {\it for some positive integer} $k$ {\it with a basis} $\{\vs_j\mid j\in\ol{-1,p^k-2}\}$ {\it satisfying}
$$\vs_{j_1}\circ \vs_{j_2}=(^{\;j_1+j_2+1}_{\;\;\;\;\;\;j_2})\vs_{j_1+j_2}+\dlt_{j_1,-1}\dlt_{j_2,0}a\vs_{p^k-2}+\dlt_{j_1,-1}\dlt_{j_2,-1}b\vs_{p^k-2} \eqno(2.2)$$
{\it for} $j_1,j_2\in\ol{-1,p^k-1}$, {\it where} $a,b\in\Bbb{F}$ {\it are constants}. 
\psp

Next, we consider infinite-dimensional simple Novikov algebra with $\mbox{char}\:\Bbb{F}=0$. Let $\Dlt$ be an additive subgroup of $\Bbb{F}$ and let $\Bbb{F}_1$ be an extension field of $\Bbb{F}$. Take a map $f:\Dlt\times\Dlt\rta \Bbb{F}_1^{\times}=\Bbb{F}_1\setminus \{0\}$ such that
$$f(\al,\be)=f(\be,\al),\;\;\;f(\al,\be)f(\al+\be,\gm)=f(\al,\be+\gm)f(\be,\gm)\eqno(2.3)$$
for $\al,\be,\gm\in\Dlt$. Take $J$ to be the additive semi-group $\{0\}$ or $\Bbb{N}$. Let ${\cal A}(\Dlt,f,J)$ be a vector space over $\Bbb{F}_1$ with a basis 
$$\{u_{\al,i}\mid \al\in\Dlt,\;i\in\Bbb{N}\}.\eqno(2.4)$$
We define an algebraic operation ``$\cdot$'' on ${\cal A}(\Dlt,f,J)$ by
$$u_{\al,i}\cdot u_{\be,j}=f(\al,\be)u_{\al+\be,i+j}\qquad\for\;\;\al,\be\in\Dlt,\;i,j\in J.\eqno(2.5)$$
Then $({\cal A}(\Dlt,f,J),\cdot)$ forms a commutative associative algebra over $\Bbb{F}_1$, and also over $\Bbb{F}$. We define a derivation $\ptl$ of ${\cal A}(\Dlt,f,J)$ over $\Bbb{F}_1$ by
$$\ptl(u_{\al,i})=\al u_{\al,i}+iu_{\al,i-1}\qquad\for\;\;\al\in\Dlt,\;i\in J,\eqno(2.6)$$
where we treat
$$u_{\be,j}=0\qquad\mbox{if}\;\;(\be,j)\not\in\Dlt\times J.\eqno(2.7)$$
Note that $\ptl$ is also a derivation of ${\cal A}(\Dlt,f,J)$ over $\Bbb{F}$.
\psp

{\bf Theorem 2.3 (Xu, [X2])}. {\it For any element} $\xi\in ({\cal A}(\Dlt,f,J)$, {\it the Novikov algebra} $({\cal A}(\Dlt,f,J),\circ_{\xi})$ {\it (cf. (2.1)) is simple}.
\psp

Osborn [O3] gave a classification of infinite-dimensional simple Novikov algebras with an element $e$ such that $e\circ e\in\Bbb{F}e$,  assuming the existence of generalized-eigenspace decomposition with respect to its left multiplication operator. There are four fundamantal mistakes in his classification. The first is using of Proposition 2.6 (d) in [O1] with $\be\neq 0$, which was misproved. The second is that the eigenspace $A_0$ in Lemma 2.12 of [O3] does not form a field when $b=0$ with respect to the Novikov algebraic operation. The third is that $A_0$ may not be a perfect field when $b\neq 0$. The fourth is that the author forgot the case $b=0$ and  $\Dlt=\{0\}$ in Lemma 2.8. In addition to these four mistakes, there are gaps in the arguments of classification in [O3]. It seems that one can not draw any conclusions of classification based on the arguments in [O3]. 

A linear transformation $T$ of a vector space $V$ is called {\it locally finite}\index{locally finite} if the subspace
$$\sum_{m=0}^{\infty}\Bbb{F}T^m(v)\;\;\mbox{is finite-dimensional for any}\;\;v\in V.\eqno(2.8)$$
An element $u$ of a Novikov algebra ${\cal N}$ is called {\it left locally finite} if its left multiplication operator $L_u$ is locally finite.
\psp

I have re-established the following classification theorem:
\psp

{\bf Theorem 2.4 (Xu, [X6])}  {\it Let} $({\cal N},\circ)$ {\it be an infinite-dimensional simple Novikov algebra over an algebraically field} $\Bbb{F}$ {\it with characteristic} $0$. {\it Suppose that} ${\cal N}$ {\it contains a left locally finite element} $e$ {\it whose right multiplication} $R_e$ {\it  is a constant map and left multiplication is surjective if} $R_e=0$. {\it Then there exist an additive subgroup} $\Dlt$ {\it of} $\Bbb{F}$, {\it an extension field} $\Bbb{F}_1$ {\it of} $\Bbb{F}$, {\it a map} $f:\Dlt\times\Dlt\rta \Bbb{F}_1^{\times}$ {\it satisfying (2.3) and} $\xi\in \Bbb{F}$ {\it such that the algebra} $({\cal N},\circ)$ {\it is isomorphic to} $({\cal A}(\Dlt,f,J),\circ_{\xi})$.
\psp

Now we consider representations of a Novikov algebra. A module $M$ of a Novikov algebras $({\cal N},\circ)$  is a vector space with two linear maps
$${\cal N}\times M\rta M: (u,w)\mapsto u\circ w,\;\;M\times {\cal N}\rta M: (w,u)\mapsto w\circ u\eqno(2.9)$$
such that (1.11) and (1.12) hold when one of the elements in $\{u,v,w\}$ is in $M$ and the other two are in ${\cal N}$. The commutator algebra associated with a finite-dimensional simple Novikov algebra over  an algebraically closed field $\Bbb{F}$ with characteristic $p>2$ is a rank-one simple Lie algebra of Witt type (cf. [O1]), whose finite-dimensional irreducible modules has not been completely classified yet. However, we have the following complete result.
\psp

{\bf Theorem 2.5 (Xu, [X1])}.  {\it Suppose that} $\Bbb{F}$ {\it is an algebraically closed field with characteristic} $p>2.$ {\it If} $M$ {\it is a finite-dimensional irreducible module of a finite-dimensional simple Novikov algebra} $({\cal N},\circ)$ {\it in the presentation (2.2),  then there exists a constant} $\lmd$ {\it and a basis} $\{v_j\mid j\in\ol{-1,p^k-2}\}$ {\it of} $M$ {\it such  that}
$$\vs_{j_1}\circ v_{j_2}=(^{\;j_1+j_2+1}_{\;\;\;\;\;\:j_2})v_{j_1+j_2}+(^{\;j_1+j_2+2}_{\;\;\;\:j_1+1})\lmd v_{j_1+j_2+1}+\dlt_{j_1,-1}\dlt_{j_2,0}av_{p^k-2}\eqno(2.10)$$
$$v_{j_2}\circ \vs_{j_1}=(^{\;j_1+j_2+1}_{\;\;\;\;\;j_1})v_{j_1+j_2}+\dlt_{j_2,-1}\dlt_{j_1,0}av_{p^k-2}+\dlt_{j_2,-1}\dlt_{j_1,-1}bv_{p^k-2} \eqno(2.11)$$
{\it for} $j_1,j_2\in\ol{-1,p^k-2}$. {\it Moreover,} $\lmd\neq 0$ {\it if} $a\neq 0$. {\it Conversely, (2.10) and (2.11) define an irreducible module for any} $0\neq \lmd \in \Bbb{F}$ {\it and for} $\lmd=0$ {\it if} $a=0$.
\psp

Let us go to the case of $\mbox{char}\:\Bbb{F}=0$. First we shall give a construction of irreducible modules. Take $J$ to be the additive semi-group $\{0\}$ or $\Bbb{N}$. Let ${\cal A}$ be a vector space with a basis 
$$\{u_{\al,i}\mid \al\in\Bbb{F},i\in J\}.\eqno(2.12)$$ 
 Define the operation ``$\cdot$'' on ${\cal  A}$ by
$$u_{\al,i}\cdot u_{\be,j}=u_{\al+\be,i+j}\qquad\for\;\;\al,\be\in \Bbb{F},\;i,j\in J.\eqno(2.13)$$
Then $({\cal  A},\cdot)$ forms a commutative associative algebra with the identity element $1=u_{0,0}$. We define the derivation $\ptl$ on ${\cal  A}$ by (2.6). Let $\Dlt$ be an additive subgroup of $\Bbb{F}$ such that $J+\Dlt\neq \{0\}$. Set
$${\cal N}=\sum_{\al\in\Dlt,i\in J}\Bbb{F}u_{\al,i}.\eqno(2.14)$$

For any fixed element $\xi\in {\cal  N}$, we define the operation ``$\circ$'' on ${\cal  A}$ by
$$u\circ v=u\cdot \ptl(v)+\xi\cdot u\cdot v\qquad\;\;\for\;\;u,v\in {\cal  A}.\eqno(2.15)$$
By Theorem 2.3, $({\cal A},\circ)$ forms a simple Novikov algebra and $({\cal N},\circ)$ forms a simple subalgebra of $({\cal A},\circ)$. For $\lmd\in\Bbb{F}$, we set
$$M(\lmd)=\sum_{\al\in\Dlt,i\in J}\Bbb{F}u_{\al+\lmd,i}\eqno(2.16)$$
Expression (2.15) shows 
$${\cal N}\circ M(\lmd),\;M(\lmd)\circ {\cal N}\subset M(\lmd).\eqno(2.17)$$
Thus $M(\lmd)$ forms an ${\cal N}$-module. By the proof of Theorem 2.9 in [X2], we have
\psp

{\bf Theorem 2.6 (Xu, [X6])}. {\it The} ${\cal N}$-{\it module} $M(\lmd)$ {\it is irreducible}.
\psp

A natural question is to what extent the modules $\{M(\lmd)\mid \lmd\in\Bbb{F}\}$ cover the irreducible modules of ${\cal N}$. Up to this point, we are not be able to answer this for a general element $\xi\in{\cal N}$. For an ${\cal N}$-module $M$, we define the left action $L_M(u_{0,0})$ of $u_{0,0}$ by 
$$ L_M(u_{0,0})(w)=u_{0,0}\circ w\qquad\for\;\;w\in M.\eqno(2.18)$$

\psp

{\bf Theorem 2.7 (Xu, [X6])}. {\it If} $\xi\in\Bbb{F}$, {\it then any irreducible} ${\cal N}$-{\it module} $M$ {\it with locally finite} $L_M(u_{0,0})$ {\it is isomorphic to} $M(\lmd)$ {\it for some} $\lmd\in\Bbb{F}$, {\it when} $\Bbb{F}$ {\it is algebraically closed.}

\section{Constructions of Gel'fand-Dorfman Bialgebras}

Recall the compatibility condition (1.17). Let $({\cal N},\circ)$ be a Novikov algebra. Define
$$[u,v]^-=u\circ v-v\circ u\qquad\for\;\;u,v\in{\cal N}.\eqno(3.1)$$

{\bf Theorem 3.1 (Gel'fand and Dorfman, [GDo])}.  {\it The triple} $({\cal N},[\cdot,\cdot]^-,\circ)$ {\it forms a Gel'fand-Dorfman bialgebra}.
\psp

Let $({\cal A},\cdot)$ be a commutative associative algebra. Denote by $\mbox{Der}\:{\cal A}$ the Lie algebras of all the derivations of ${\cal A}$.
Set
$${\cal N}=\mbox{Der}\:{\cal A}\oplus {\cal A}.\eqno(3.2)$$
We define two algebraic operation $[\cdot,\cdot]$ and $\circ$ on ${\cal N}$ by
$$[d_1+\xi_1,d_2+\xi_2]=[d_1,d_2]+d_1(\xi_2)-d_2(\xi_1),\eqno(3.3)$$
$$(d_1+\xi_1)\circ(d_2+\xi_2)=\xi_2(d_1+\xi_1)\eqno(3.4)$$
for $d_1,d_2\in \mbox{Der}\:{\cal A}$ and $\xi_1,\xi_2\in{\cal A}$.
\psp

{\bf Theorem 3.2 (Xu, [X4])}. {\it The triple} $({\cal N},[\cdot,\cdot],\circ)$ {\it forms a Gel'fand-Dorfman  bialgebra}.
\psp

The above construction is extracted from the simple Lie algebras of Witt type.

Our second construction is related to the following concept. A {\it Lie-Poisson algebra} is a vector space ${\cal A}$  with two algebraic operations ``$\cdot$'' and $[\cdot,\cdot]$ such that $({\cal A},\cdot)$ forms a commutative associative algebra, $({\cal A},[\cdot,\cdot])$ forms a Lie algebra and the following compatibility condition is satisfied:
$$[u,v\cdot w]=[u,v]\cdot w+v\cdot [u,w]\qquad\mbox{for}\;\;u,v,w\in{\cal A}.\eqno(3.5)$$

 Let $({\cal A},\cdot,[\cdot,\cdot])$ be a Lie-Poisson algebra and let $\ptl$ be a derivation of the algebra $({\cal A},\cdot)$ such that
$$\ptl[u,v]=[\ptl(u),v]+[u,\ptl(v)]+\xi[u,v]\qquad\for\;\;u,v\in{\cal A},\eqno(3.6)$$
where $\xi\in \Bbb{F}$ is a constant.
Now we define another algebraic operation on $\circ$ on ${\cal A}$ by
$$u\circ v=u\ptl(v)+\xi uv\qquad\mbox{for}\;\;u\in {\cal A},\;v\in {\cal A}.\eqno(3.7)$$

{\bf Theorem 3.3 (Xu, [X4])}. {\it The triple} $({\cal A},[\cdot,\cdot],\circ)$ {\it forms a Gel'fand-Dorfman  bialgebra}.

The above construction is related to the Lie algebras of Hamiltonian type and Contact type.

Let $({\cal A},\cdot)$ be a commutative associative algebra and let $\ptl_1,\ptl_2$ be mutually commutative derivations  of $({\cal A},\cdot)$. Define algebraic operations on ${\cal A}$ by
$$[u,v]=\ptl_1(u)\ptl_2(v)-\ptl_2(u)\ptl_1(v)+u\ptl_2(v)-\ptl_2(u)v,\;\;u\circ v=u\ptl_2(v)\eqno(3.8)$$
for $u,v\in{\cal A}$. 
\psp

{\bf Theorem 3.4 (Xu, [X4])}. {\it The triple} $({\cal A},[\cdot,\cdot],\circ)$ {\it forms a Gel'fand-Dorfman bialgebra}.
\psp

Let $({\cal A},\cdot),\;\ptl_1$ and $\ptl_2$ be the same as in the above. For any constant $b\in\Bbb{F}$, we define 
$$[u,v]=\ptl_1(u)\ptl_2(v)-\ptl_2(u)\ptl_1(v)+b(u\ptl_2(v)-\ptl_2(u)v),\;\;u\circ v=u\ptl_1(v)+buv\eqno(3.9)$$
for $u,v\in{\cal A}$. 
\psp

{\bf Theorem 3.5 (Xu, [X4])}. {\it The triple} $({\cal A},[\cdot,\cdot],\circ)$ {\it forms a Gel'fand-Dorfman bialgebra}.

\section{Classifications of the Bialgebras}

If $({\cal N},[\cdot,\cdot],\circ)$ is a Gel'fand-Dorfman bialgebra, then we say that $({\cal N},[\cdot,\cdot])$ is a {\it Lie algebra over the Novikov algebra} $({\cal N},\circ)$. Similarly, $({\cal N},\circ)$ is a {\it Novikov algebra over the Lie algebra} $({\cal N},[\cdot,\cdot])$. 

 Let $\Bbb{F}$ be an algebraically closed field with $\mbox{char}\:\Bbb{F}=0$. Take $J$ to be the additive semi-group $\{0\}$ or $\Bbb{N}$ and take an additive subgroup $\Dlt$ of $\Bbb{F}$.  Let ${\cal A}(\Dlt,J)$ be a vector space with a basis  
$$\{u_{\al,i}\mid \al\in\Dlt,\;i\in J\}.\eqno(4.1)$$ 
For any constant $b\in\Bbb{F}$, we define an algebraic operation $\circ_b$ on 
${\cal A}(\Dlt,J)$ by
$$u\circ_b v=u\cdot \ptl(v)+b\cdot u\cdot v\qquad\;\;\for\;\;u,v\in {\cal  A}(\Dlt,J).\eqno(4.2)$$
By Theorem 2.3, $({\cal A}(\Dlt,J),\circ_b)$ forms a simple Novikov algebra.
\psp

{\bf Theorem 4.1 (Osborn and Zel'manov, [OZ])}. {\it If} $\Dlt=\Bbb{Z},\;J=\{0\},\; b\not\in\Bbb{Z}$ {\it or} $\Dlt=0,\;J=\Bbb{N}$, {\it then any nontrivial Lie algebra over the Novikov algebra} $({\cal A}(\Dlt,J),\circ_b)$ {\it is isomorphic to} $({\cal A}(\Dlt,J),[\cdot,\cdot]^-)$, {\it where}
$$[u_{\al,i},u_{\be,j}]^-=(\be-\al)u_{\al+\be,i+j}+ (j-i)u_{\al+\be,i+j-1}\eqno(4.3)$$
{\it for} $\al,\be\in\Dlt$ {\it and} $i,j\in J$. 
\pse

{\bf Theorem 4.2 (Xu, [X4])}. {\it If} $b\not\in\Dlt$, {\it then any Lie algebra over the  Novikov algebra} $({\cal A}(\Dlt,J),\circ_b)$ {\it has the Lie bracket}
\begin{eqnarray*}[u_{\al,i},u_{\be,j}]&=&[(\al+b)\vf(\be)-(\be+b)\vf(\al)]u_{\al+\be,i+j}\\& &+[i(\vf(\be)-\lmd(\be+b))+j(\lmd(\al+b)-\vf(\al))]u_{\al+\be,i+j-1}\hspace{3.4cm}(4.4)\end{eqnarray*}
{\it for} $\al,\be\in\Delta,\;i,j\in J$, {\it where} $\vf:\Delta\rightarrow \Bbb{F}$ {\it is an additive group homomorphism and} $\lmd\in\Bbb{F}$ {\it is a constant}.
\psp

When $J=\{0\}$, (4.4) gives a Block algebra (cf. [B]). Next we assume $J=\{0\}$ and redenote $u_{\al,0}$ by $u_{\al}$ for $\al\in\Dlt$.
\psp

{\bf Theorem 4.3 (Xu, [X4])} {\it Any Lie algebra over the Novikov algebra} $({\cal A}(\Dlt,\{0\}),\circ_0)$ {\it has the Lie bracket}
$$[u_{\al},u_{\be}]=(\phi(\al,\be)+a(\be-\al))u_{\al+\be}\qquad\mbox{\it for}\;\;\al,\be\in \Dlt,\eqno(4.5)$$
{\it where} $a$ {\it is a constant and} $\phi(\cdot,\cdot):\Dlt\times\Dlt\rightarrow \Bbb{F}$ {\it is a skew-symmetric map such that there exists a symmetric map} $S_0(\cdot,\cdot,\cdot):\Dlt\times\Dlt\times\Dlt\rightarrow \Bbb{F}$ {\it satisfying} 
$$\phi(\be+\gm,\al)=\phi(\gm,\al)+\phi(\be,\al)+\al S_0(\al,\be,\gm),\eqno(4.6)$$
$$(\gm\phi(\al,\be)+\al\phi(\be,\gm)+\be\phi(\gm,\al))(S_0(\al,\be,\gm)-a)=0.\eqno(4.7)$$
{\it for} $\al,\be,\gm\in\Dlt.$ {\it In particular, the following} $\phi$ {\it satisfies our condition: when} $a=0$, $\phi$ {\it is any skew-symmetric} $\Bbb{Z}$-{\it bilinear form; when} $a\neq 0$,
$$\phi(\al,\be)=\al\vf_0(\be)-\be\vf_0(\al)\qquad\mbox{\it for}\;\;\al,\be\in\Dlt,\eqno(4.8)$$
{\it where} $\vf_0:\Dlt\rightarrow \Bbb{F}$ {\it is an additive group homomorphism}. 
\psp

{\bf Theorem 4.4 (Xu, [X4])} {\it Suppose} $0\neq b\in\Dlt$. {\it
 Any Lie algebra over the  Novikov algebra} $({\cal A}(\Dlt,\{0\}),\circ_b)$ {\it has the Lie bracket}
$$[u_{\al},u_{\be}]=\theta(\al,\be)u_{\al+\be+b}+((\al+b)\vf(\be)-(\be+b)\vf(\al))u_{\al+\be}\eqno(4.9)$$
{\it for} $\al,\be\in\Dlt$, {\it where} $\vf:\Dlt\rightarrow \Bbb{F}$ {\it is additive group homomorphism,} $\theta\equiv 0$ {\it if} $\vf(b)\neq 0$ {\it and} $\theta:\Dlt\times\Dlt\rightarrow\Bbb{F}$ {\it is a skew-symmetric map satisfying}
$$(\al+b)(\theta(\al+\gm,\be)-\theta(\gm,\be)-\theta(\al,\be))=(\be+b)(\theta(\be+\gm,\al)-\theta(\gm,\al)-\theta(\be,\al))\eqno(4.10)$$
{\it and}
$$\theta(\al,\be)\theta(\al+\be+b,\gm)+\theta(\be,\gm)\theta(\be+\gm+b,\al)+\theta(\gm,\al)\theta(\gm+\al+b,\be)=0\eqno(4.11)$$
{\it if} $\vf(b)=0$. {\it In particualr, we can take} $\vf:\Dlt\rightarrow \Bbb{F}$ {\it to be any  additive group homomorphism}  {\it such that} $\vf(b)=0$ {\it and take} $\theta(\cdot,\cdot)${\it to be  a skew-symmetric} $\Bbb{Z}${\it -bilinear map such that} $b\in\mbox{Rad}_{\theta}$ {\it or}
$$\theta(\al,\be)=\vf_1(\al)\vf_2(\be)-\vf_1(\be)\vf_2(\al)\qquad{\it for}\;\;\al,\be\in\Delta,\eqno(4.12)$$
{\it where} $\vf_1,\vf_2:\Delta\rightarrow \Bbb{F}$ {\it are additive group homorphisms such that} $\vf_1(b)=0$ {\it and} $\vf_2(b)\neq 0$.
\psp

{\bf Theorem 4.5 (Xu, [X4])}. {\it If} $\circ$ {\it is a Novikov algebra operation on the space}  ${\cal A}(\Delta,\{0\})$ {\it such that}
$$u_{\al}\circ u_{\be}-u_{\be}\circ u_{\al}=(\be-\al)u_{\al+\be}\qquad\mbox{\it for}\;\;\al,\be\in\Dlt,\eqno(4.13)$$
{\it then there exists an element} $\xi\in  {\cal A}(\Delta,\{0\})$ {\it such that}
$$u_{\al}\circ u_{\be}=(\be+\xi)u_{\al+\be}\qquad\mbox{for}\;\;\al,\be\in\Dlt.\eqno(4.14)$$

\section{Application to Cubic and Quartic Conformal Algebras}

 For two vector spaces $V_1$ and $V_2$, we denote by $LM(V_1,V_2)$ the space of linear maps from $V_1$ to $V_2$. We shall also use the following operator of taking residue:
$$\rd_z(z^n)=\dlt_{n,-1}\qquad\for\;\;n\in \Bbb{Z}.\eqno(5.1)$$
Furthermore, all the binomials are assumed to be expanded in the nonnegative powers of the second variable. For example
,
$${1\over z-x}={1\over z(1-x/z)}=\sum_{j=0}^{\infty}z^{-1}\left({x\over z}\right)^j=\sum_{j=0}^{\infty}z^{-j-1}x^j.\eqno(5.2)$$
In particular, the above equation implies
$$\rd_x{1\over z-x}(\sum_{j\in\Bbb{Z}}\xi_jz^j)=\sum_{j=1}^{\infty}\xi_{-j}z^{-j}.\eqno(5.3)$$
So the operator $\rd_x(1/( z-x))$(---) is taking the part of negative powers in a formal series and changing the variable $x$ to $z$.
\psp

 A {\it conformal algebra} $R$ is  an $\Bbb{F}[\ptl]$-module equipped with a  linear map $Y^+(\cdot,z):\;R\rightarrow LM(R,R[z^{-1}]z^{-1})$ satisfying:
$$Y^+(\ptl u,z)={d\over dz}Y^+(u,z),\;\;Y^+(u,z)v=\rd_x{e^{x\ptl}Y^+(v,-x)u\over z-x},\eqno(5.4)$$
$$Y^+(u,z_1)Y^+(v,z_2)-Y^+(v,z_2)Y^+(u,z_1)=\rd_x{Y^+(Y^+(u,z_1-x)v,x)\over z_2-x} \eqno(5.5)$$
for $u,v\in R$. We denote by $(R,\ptl,Y^+(\cdot,z))$ a conformal algebra.
\psp

The above definition is the equivalent generating-function form to that given in [K], where the author used the component formula with $Y^+(u,z)=\sum_{n=0}^{\infty}u_{(n)}z^{-n-1}/n!$. The connection between between the Lie algebra with one-variable structure in (1.3) and conformal algebra is that $R=\Bbb{F}[\ptl]\otimes_{\Bbb{F}}V$ and
$$Y^+(u,z)v=\sum_{i=0}^m\sum_{j=0}^n(-1)^ii!\ptl^jw_{ij}z^{-i-1}.\eqno(5.6)$$
Conformal algebras are equivalent to linear Hamiltonian operators (cf. [X5]).

 Suppose that $(R,\ptl,Y^+(\cdot,z))$ is a conformal algebra that is a free $\Bbb{F}[\ptl]$-module over a subspace $V$, namely
$$R=\Bbb{F}[\ptl]V\;\;(\cong \Bbb{F}[\ptl]\otimes_{\Bbb{F}}V).\eqno(5.7)$$
Let $m$ be a positive integer. The algebra $R$ is called  {\it of degree} $m$ if for any $u,v\in V$,
$$Y^+(u,z)v=\sum_{0<j;i+j\leq m}\ptl^iw_{i,j}z^{-j}\qquad\mbox{with}\;\;w_{i,j}\in V,\eqno(5.8)$$
and $w_{m-j,j}\neq 0$ for some $u,v\in V$ and $j\in\ol{1,m}$. A {\it quadratic conformal algebra} is a conformal algebra of degree 2, a {\it cubic conformal algebra} is a conformal algebra of degree 3 and a {\it quartic conformal algebra} is a conformal algebra of degree 4. In [X4], I have proved that a quadratic conformal algebra is equivalent to a Gel'fand-Dorfman bialgebra. Below I will give the construction of quadratic conformal algebra from Gel'fand-Dorfman bialgebras.

Let $({\cal N},[\cdot,\cdot],\circ)$ be a Gel'fand Dorfman bialgebra. Set
$$R_{\cal N}=\Bbb{F}[\ptl]\otimes_{\Bbb{F}}{\cal N}.\eqno(5.9)$$
So $R_{\cal N}$ is a free $\Bbb{F}[\ptl]$-module generated by ${\cal N}$. For convenience, we identify ${\cal N}$ with $1\otimes {\cal N}$. We define a  linear map $Y^+(\cdot,z):\;R\rightarrow LM(R,R[z^{-1}]z^{-1})$ by
$$Y^+(\ptl^mu,z)\ptl^nv=\sum_{j=0}^n(-1)^j(^n_j)\left({d\over dz}\right)^{m+j}
\ptl^{n-j}(([v,u]+\ptl(v\circ u))z^{-1}+(u\circ v+v\circ u)z^{-2})\eqno(5.10)$$
for $u,v\in {\cal N}$ and $m,n\in\Bbb{N}$.
\psp

{\bf Theorem 5.1 (Xu, [X4])}. {\it The triple} $(R_{\cal N},\ptl,Y^+(\cdot,z))$ {\it forms a quadratic conformal algebra}.
\psp

Next I will use the above theorem to construct  simple cubic conformal algebras and quartic conformal algebras. A conformal algebra $(R,\ptl,Y^+(\cdot,z))$ is called {\it simple} if there does not exist a nonzero proper subspace ${\cal I}$ of $R$ such that
$$Y^+(u,z)({\cal I})\subset {\cal I}[z^{-1}]\qquad\for\;\; u\in R.\eqno(5.11)$$

Take $J$ to be the additive semi-group $\{0\}$ or $\Bbb{N}$. Let $\G$ be an additive subgroup of $\Bbb{F}^{\:2}$ such that
$$(J,0)+\G\not\subset (0,\Bbb{F}),\qquad \G\not\subset (\Bbb{F},0).\eqno(5.12)$$
Let ${\cal A}$ be the semi-group algebra of $\G\times J$ with the canonical basis $\{u_{\al,i}\mid\al\in\Dlt,\;i\in J\}$, that is,
$$u_{\al,i}\cdot u_{\be,j}=u_{\al+\be,i+j}\qquad\for\;\; \al,\be\in\G,\;i,j\in J.\eqno(5.13)$$
Define two derivations $\ptl_1$ and $\ptl_2$ of ${\cal A}$ by
$$\ptl_1(u_{\al,i})=\al_1u_{\al,i}+iu_{\al,i-1},\;\;\;\ptl_2(u_{\al,i})=\al_2u_{\al,i}\eqno(5.14)$$
for $\al=(\al_1,\al_2)\in\G$ and $i\in J$. 

First we define
$$[u,v]_1=\ptl_1(u)\ptl_2(v)-\ptl_2(u)\ptl_1(v)+u\ptl_2(v)-\ptl_2(u)v,\;\;u\circ_1 v=u\ptl_2(v)\eqno(5.15)$$
Then $({\cal A},[\cdot,\cdot]_1,\circ_1)$ forms a Gel'fand-Dorfman bialgebra by Theorem 3.4. On
$R_{\cal A}$, the structure map (5.10) is determined by
\begin{eqnarray*}Y^+_1(u_{\al,i},z)u_{\be,j}&=&[((1+\al_1)\be_2-\al_2(1+\be_1))u_{\al+\be,i+j}+(i\be_2-j\al_2)u_{\al+\be,i+j-1}\\& &+\be_2\ptl u_{\al+\be,i+j}]z^{-1}+(\al_2+\be_2)u_{\al+\be,i+j}z^{-2}.\hspace{4.5cm}(5.16)\end{eqnarray*}
for $\al=(\al_1,\al_2),\;\be=(\be_1,\be_2)\in\G$ and $i,j\in J$. When $\be_2=-\al_2$, (5.16) becomes
$$Y^+_1(u_{\al,i},z)u_{\be,j}=\be_2[(\al_1+\be_1+2+\ptl)u_{\al+\be,i+j}+(i+j)u_{\al+\be,i+j-1}]z^{-1}.\eqno(5.17)$$
Set
$$R_1=\sum_{\al=(\al_1,\al_2)\in\G,\:\al_2\neq 0}\Bbb{F}[\ptl]u_{\al,i}+\sum_{\be=(\be_1,0)\in\G,\;j\in J}\Bbb{F}[\ptl][(\be_1+2+\ptl)u_{\be,j}+ju_{\be,j-1}].\eqno(5.18)$$
Then $R_1$ forms a subalgebra of $(R_{\cal A},\ptl,Y^+_1(\cdot,z))$ by (5.17), that is,
$$\ptl R_1\subset R_1,\;\;Y^+_1(u,z)v\in R_1[z^{-1}]\qquad\for\;\;u,v\in R_1.\eqno(5.19)$$

{\bf Theorem 5.2 (Xu, [X7])}. {\it The conformal algebra} $(R_1,\ptl,Y^+_1(u,z))$ {\it is a simple cubic  conformal algebra}.
\psp

Let ${\cal A},\;\ptl_1$ and $\ptl_2$ be the same as in the above. Take $b\in\Bbb{F}$ to be any fixed constant such that 
$$(\Bbb{F},2b)\bigcap\G\neq \emptyset,\qquad \G\not\subset (\Bbb{F},0).\eqno(5.20)$$
 We define
$$[u,v]_2=\ptl_2(u)\ptl_1(v)-\ptl_1(u)\ptl_2(v)+b(u\ptl_1(v)-\ptl_1(u)v,\;\;u\circ_2 v=u\ptl_2(v)+buv\eqno(5.21)$$
for $u,v\in{\cal A}$.

Then $({\cal A},[\cdot,\cdot]_2,\circ_2)$ forms a Gel'fand-Dorfman bialgebra by Theorem 3.5. On $R_{\cal A}$, the structure map (5.10) is determined by
\begin{eqnarray*}Y^+_2(u_{\al,i},z)u_{\be,j}&=&[(\al_1(\be_2+b)-(\al_2+b)\be_1))u_{\al+\be,i+j}+(i(\be_2+b)-j(\al_2+b))u_{\al+\be,i+j-1}\\& &+(\be_2+b)\ptl u_{\al+\be,i+j}]z^{-1}+(\al_2+\be_2+2b)u_{\al+\be,i+j}z^{-2}.\hspace{2.4cm}(5.22)\end{eqnarray*}
for $\al=(\al_1,\al_2),\;\be=(\be_1,\be_2)\in\G$ and $i,j\in J$. When $\al_2=2b-\be_2$, (5.22) becomes
$$Y^+_2(u_{\al,i},z)u_{\be,j}=(\be_2+b)[(\al_1+\be_1+\ptl)u_{\al+\be,i+j}+(i+j)u_{\al+\be,i+j-1}]z^{-1}.\eqno(5.23)$$

Set
$$R_2=\sum_{\al=(\al_1,\al_2)\in\G,\:\al_2\neq-2b}\Bbb{F}[\ptl]u_{\al,i}+\sum_{\be=(\be_1,-2b)\in\G,\;j\in J}\Bbb{F}[\ptl][(\be_1+\ptl)u_{\be,j}+ju_{\be,j-1}].\eqno(5.24)$$
Then $R_2$ forms a subalgebra of $(R_{\cal A},\ptl,Y^+_2(\cdot,z))$ by (5.23), that is,
$$\ptl R_2\subset R_2,\;\;Y^+_2(u,z)v\in R_2[z^{-1}]\qquad\for\;\;u,v\in R_2.\eqno(5.25)$$

{\bf Theorem 5.3 (Xu, [X7])}. {\it The conformal algebra} $(R_2,\ptl,Y^+_2(u,z))$ {\it is a simple  cubic conformal algebra if} $b=0$ {\it and is a simple quartic conformal algebra if} $b\neq 0$.
\vspace{1cm}

\noindent{\Large \bf References}

\hspace{0.5cm}

\begin{description}

\item[{[A]}] L. Ausland, Simple transitive groups of affine motions, {\it Amer. J. Math.} {\bf 99} (1977), 809-826.

\item[{[BN]}] A. A. Balinskii and S. P. Novikov, Poisson brackets of hydrodynamic type, Frobenius algebras and Lie algebras, {\it Soviet Math. Dokl.} Vol. {\bf 32} (1985), No. {\bf 1}, 228-231.

\item[{[B]}] R. Block, On torsion-free abelian groups and Lie algebras, {\it Proc. Amer. Math. Soc.} {\bf 9} (1958), 613-620.

\item[{[F]}] V. T. Filipov, A class of simple nonassociative algebras, {\it Mat. Zametki} {\bf 45} (1989), 101-105.

\item[{[FG]}] D. Fried and W. Goldman, Three dimensional affine crystallographic groups, {\it Adv. Math}. {\bf 47} (1983), 1-49.

\item[{[GDi]}] I. M. Gel'fand and L. A. Dikii, Asymptotic behaviour of the resolvent of Sturm-Liouville equations and the algebra of the Korteweg-de Vries equations, {\it Russian Math. Surveys} {\bf 30:5} (1975), 77-113.

\item[{[GDo]}] I. M. Gel'fand and I. Ya. Dorfman, Hamiltonian operators and algebraic structures related to them, {\it Funkts. Anal. Prilozhen}  {\bf 13} (1979), 13-30.

\item[{[K]}] V. G. Kac, {\it Vertex algebras for beginners}, University lectures series, Vol {\bf 10}, AMS. Providence RI, 1996.

\item[{[O1]}]
J. Marshall Osborn, Novikov algebras, {\it Nova J. Algebra} \& {\it Geom.} {\bf 1} (1992), 1-14.

\item[{[O2]}]
---, Simple Novikov algebras with an idempotent, {\it Commun. Algebra} {\bf 20} (1992), No. 9, 2729-2753.

\item[{[O3]}]
---, Infinite dimensional Novikov algebras of characteristic 0, {\it J. Algebra} {\bf 167} (1994), 146-167.

\item[{[O4]}]---, Modules for Novikov algebras, {\it Proceeding of the II International Congress on Algebra, Barnaul, 1991.}

\item[{[O5]}]
---, Modules for Novikov algebras of characteristic 0, {\it Commun. Algebra} {\bf 23} (1995), 3627-3640.

\item[{[OZ]}] J. Marshall Osborn and E. Zel'manov, Nonassociative algebras related to Hamiltonian operators in the formal calculus of variations, {\it J. Pure. Appl. Algebra} {\bf 101} (1995), 335-352.

\item[{[X1]}] X. Xu, On simple Novikov algebras and their irreducible modules, {\it J. Algebra} {\bf 185} (1996), 905-934.

\item[{[X2]}]---, Novikov-Poisson algebras, {\it J. Algebra} {\bf 190} (1997), 253-279.

\item[{[X3]}] ---, Variational calculus of supervariables and related algebraic structures, {\it J. Algebra} {\bf 223} (2000), 396-437.

\item[{[X4]}] ---,  Quadratic conformal superalgebras, {\it J. Algebra} {\bf 231} (2000), 1-38.

\item[{[X5]}]---, Equivalence of conformal superalgebras to Hamiltonian superoperators, to appear in {\it Algebra Colloquium} {\bf 8} (2001), 63-92.

\item[{[X6]}]---, On classification of simple Novikov algebras and their irreducible modules over a field with characteristic 0,  ArXiv:Math.QA/0008072, 2000.

\item[{[X7]}]---,  {\it Algebraic Theory of Hamiltonian Superoperators}, a monograph, initial preprint.

\item[{[Z]}]E. I. Zel'manov, On a class of local translation invariant Lie algebras, {\it Soviet Math. Dokl.} Vol {\bf 35} (1987), No. {\bf 1}, 216-218.

\end{description}

\noindent Department of Mathematics, The Hong Kong University of Science and Technology, center, Clear Water Bay, Kowloon, Hong Kong

\end{document}